\documentclass[12pt]{article} \input{diagrams}
\begin{document}\newtheorem{defn}{Definition}[section]
\newtheorem{prop}{Proposition}[section]
\newtheorem{lemma}{Lemma}[section]
\newtheorem{thm}{Theorem}[section]

\newcommand{\cz}{{C^{\infty}}}
\newcommand{\Ci}{C^{\infty}({\bf R}^l ) }
\newcommand{\ci}{C^{\infty}}
\newcommand{\ra}{\to}
\newcommand{\lra}{\to}
\title{Categorical distribution theory; heat equation}
\author{A.  Kock \and G.E.  Reyes}
\date{}
  \maketitle

\newcommand{\Cil}{C^{\infty}({\bf R}^l)}

\section*{Introduction} The simplest notion by which a theory of function
spaces may be 
formulated is  that of cartesian closed categories.
To realize this concretely for spaces of  smooth (= $C^{\infty}$) functions,
several notions of diffeological spaces and convenient
vector spaces have been developed,  
besides  the whole body of topos theory. Topos theory 
in
particular provides for toposes containing the category of smooth
manifolds as full subcategory.  In fact, Grothendieck's ``Smooth
Topos'' is closely related to the category of diffeological spaces.
The special features of Convenient Vector Spaces were utilized by the
present authors back in the 1980's for a more elaborate topos , cf.\
\cite{CVSE}, \cite{CA}.  The topos there (Dubuc's ``Cahiers Topos'') in
fact accomodates synthetic
differential geometry.  The present work is a continuation of our
work from the 80's, and is motivated by the desire to have a synthetic theory
of some of the fundamental partial differential equations, like the
heat equation.  This forced us to sort out how distribution theory (in
the sense of L.\ Schwartz) relates to convenient vector space theory
and to the Cahiers Topos.  Note that for the heat equation,
distributions of compact support will not suffice (distributions of
compact support are easier to deal with in categorical terms, as we
did in our paper on the wave equation, \cite{WE}).

   In particular, we study smoothness with respect to time of solutions of the
   heat equation\footnote{on the unlimited line; we did not work out the
   details for higher dimensional Euclidean spaces, let alone Riemannian
   manifolds.}.  These solutions model evolution through time of a {\it
   heat distribution}.  A heat distribution is an {\it extensive}
   quantity and does not necessarily have a density {\em function},
   which is an {\it intensive} quantity.  The most important of all
   distributions, the point- or Dirac- distributions, do not.  For the
   heat equation, it is well known that the evolution through time of any
   distribution leads `instantaneously' (i.e., after any {\it positive}
   lapse of time $t>0$) to distributions that do have smooth density
   functions.  Indeed, the evolution through time of the Dirac
   distribution $\delta(0)$ is given by the map (``heat kernel'', 
   ``fundamental solution'')
   \begin{equation}K:{\bf R}_{\geq 0}\lra {\cal D}'({\bf
   R})\label{K}\end{equation} defined by cases by the classical formula
   \begin{equation} K(t)=\left\{\begin{array}{ll} e^{-x^2/4t}/\sqrt{4\pi
   t} & \mbox{if $t>0$} \\
\delta(0) & \mbox{if $t=0$}
\end{array}
\right. ;
\label{heatkernel}\end{equation} here ${\cal D}' ({\bf R}) $ denotes a
suitable space of distributions (in the sense of
\cite{Schw2}, \cite{Schw1}); notice that in the first clause we are
identifying distributions with their density functions (when such
density functions exist).

\medskip
The fundamental mathematical object given in (\ref{heatkernel})
presents a challenge to the synthetic kind of reasoning in
differential geometry, where a basic tenet is ``everything is
smooth''; therefore, definition by cases, as in (\ref{heatkernel}), has a
dubious status.  It was this challenge that motivated the present
study; more precisely, we wanted to present a model of Synthetic
Differential Geometry where the map (\ref{K}) does exist, and
satisfies the heat equation, (as well as (\ref{heatkernel})).  We shall
in fact construct such $K$ in the Cahiers topos \cite{Dubuc}.

This construction leads to some smoothness questions that have a
purely classical formulation, see Section 4 below.

\medskip One may see another lack of smoothness in
(\ref{heatkernel}), namely ``$\delta (0)$ is not smooth''; but this ``lack
of smoothness'' is completely spurious, when one firmly stays in the
space of distributions and their intrinsic ``diffeology'', in
particular avoiding to view distributions as generalized functions.  We
describe in Section 2 the distribution theory that is adequate for the
purpose.  In fact, as will be seen in Section 9, this theory is
forced on us by synthetic considerations in the Cahiers topos.

\medskip
We want to thank Henrik Stetk\ae r for useful conversations on the topic
of distributions.

\section{Diffeological spaces and convenient vector spaces}
We collect some notions and facts.  Some references are collected at the
end of the section.

A {\em diffeological space} is a
set $X$ equipped with a collection of {\em smooth plots}, a plot $p$ being
a map from (the underlying set of) an open set $U$ of some $ {\bf R}^n$
into $X$, $p: U \to X$; the collection should satisfy certain
stability properties.  These properties are best summarized by
considering the following site $\underline{mf}$: its objects are open
subsets of ${\bf R}^n$, the maps are smooth maps between such sets; a
covering is a jointly surjective family of local diffeomorphisms.  (This
site is a site of definition of the ``Smooth Topos'' of Grothendieck
et al., \cite{SGA4} p.\ 318; and is one of the first examples of what they
call a ``Gros Topos''.) Any set $X$ gives rise to a presheaf $c(X)$ on
this site, namely $c(X)(U):= Hom _{sets}(U,X)$.  A diffeological
structure on the set $X$ is a subsheaf $P$ of the presheaf $c(X)$, the
elements of $P(U)$ are called the {\em smooth} $U$-plots on $X$.  A set
theoretic map $f: X\to X'$ between diffeological spaces is called
{\em (plot-) smooth} if $f\circ p$ is a smooth plot on $X'$ whenever $p$ is a smooth
plot on $X$.

Any smooth manifold $M$ carries a canonical
diffeology, namely with $P(U)$ being the set of smooth maps $U\to M$.
We have full inclusions of categories: smooth manifolds into
diffeological spaces into the smooth topos, (= the topos of sheaves on
the site $\underline{mf}$), $$\underline{Mf} \subseteq \underline{Diff}
\subseteq
\underline{sh}(\underline{mf}).$$

If $X$ is a diffeological space, and $H\subseteq X$ a subset, there is
an induced diffeology on $H$, namely by declaring $U\to H$ to be a
smooth plot iff it is a smooth plot viewed as a map into $X$.  
In particular the non-negative
reals (=the closed half line) ${\bf R}_{\geq 0} \subseteq {\bf R}$ 
will be considered a diffeologicaal space with the diffeology induced 
by that of ${\bf R}$.

\medskip

Let ${\bf R}_{>0}$ denote the open half line of positive reals.  A
smooth function $f:{\bf R}_{>0} \to {\bf R}$ is called {\em
square-smooth} if $f(x^2 )$ is of the form $g(x)$ for a smooth
function $g: {\bf R}\to {\bf R}$ (necessarily unique).

Note that the square root function is smooth on ${\bf R}_{>0}$, but not square
smooth, since $\sqrt{x^2} = |x|$, which does not extend smoothly to the
whole line.  Note also that if $f$ is square smooth, then it extends
(uniquely) to a continuous function on the closed half line ${\bf
R}_{\geq 0}$, by putting $f(0)=g(0)$.

\begin{prop} If $f:{\bf R}_{>0}\to {\bf R}$ is square smooth, then so is $f'$.
\label{p11}\end{prop} {\bf Proof.} Let $f(x^2) = g(x)$, with $g$ smooth.  Then
clearly $g$ is an even function, so $g'(0) =0$.  Therefore, $g'(x) =
x\cdot h(x)$ for a unique smooth function $h$.  Also we note that
since $g$ is even, $g'$ is odd.  For $t>0$, we have $f(t)=
g(t^{1/2})$, and so for $t>0$, $$f'(t)= \frac{1}{2}\;
g'(t^{1/2})\cdot t^{-1/2}.$$
For $x\neq 0$, we therefore have (using $\sqrt (x^2 )= |x|$)
    $$f'(x^2 ) =
\frac{1}{2} g'(|x|)\cdot |x|^{-1},$$
but since $g'$ is odd, $g'(|x|)\cdot |x|^{-1} = g'(x)\cdot x^{-1}$.
Thus, for $x\neq 0$,
$$f'(x^2 )= \frac{1}{2} \; g'(x)\cdot x^{-1} = \frac{1}{2}\;  x\cdot
h(x)\cdot x^{-1}$$
which extends to the smooth function $1/2 \; h(x)$, defined on the whole line.
This proves that $f'$ is square smooth.

\medskip

A smooth function $f:{\bf R}_{>0}\to {\bf R}$ is called {\em Seeley-smooth}
(after \cite{Seeley}) if all the higher derivatives $f^{(k)}: {\bf R}_{>0}\to
{\bf R}$ have finite limits as $t\to 0+$, i.e.  if each $f^{(k)}$
extends to a continuous function defined on the closed half line
${\bf R}_{\geq 0}$.  A Corollary of the Proposition is then

\begin{prop} If a function $f:{\bf R}_{>0}\to {\bf R}$  is square-smooth,
it is Seeley smooth.
\label{p22}\end{prop} {\bf Proof.} We already observed that a square smooth $f
:{\bf R}_{>0} \to {\bf R}$ extends continuously to ${\bf R}_{\geq 0}$.  But
since by Proposition \ref{p11}, all higher derivatives of $f$ are also square
smooth, they  extend continuously as well.

\medskip
\begin{thm} For a function $f:{\bf R}_{\geq 0}\to {\bf R}$, the following
conditions are equivalent

1) $f$ is plot smooth;

2) the restriction of $f$ to ${\bf R}_{>0}$ is square smooth;

3) the restriction of $f$ to ${\bf R}_{>0}$ is Seeley smooth;

4) $f$ extends to a smooth function defined on the whole of ${\bf R}$.
\label{p33}\end{thm}
A function ${\bf R}_{\geq 0}\to {\bf R}$ satisfying these conditions, we will
simply call {\em smooth}.

\medskip

   {\bf Proof.} 1) implies 2), since the function
$x^2$ is one of the plots we may use; 2) implies 3), by Proposition
\ref{p22}; 3) implies 4) by Seeley's Theorem, \cite{Seeley}
(alternatively, by Borel's extension theorem for formal power series,
see e.g.  \cite{MR} p.\ 18); and clearly, the restriction of a global smooth
function to ${\bf R}_{\geq 0}$ (or to any other subset of ${\bf R}$,
with its induced diffeology) is plot smooth (= diffeologically
smooth); so 4) implies 1).  (An alternative proof of $1) \Rightarrow
4)$ follows from Whitney's theorem on even functions, cf.\ \cite{PR}.)

\medskip
Note that it follows that if $f$ is Seeley smooth, then the assumed limit
$f^{(k)}(t)$ as $t\to 0+$ is the $k$'th derivative at $0$ of (any
smooth extension of) $f$.

\medskip
    The category of diffeological spaces $\underline{Diff}$ is
cartesian closed (in fact, it is a concrete quasi-topos).  Thus, if
$X$ and $Y$ are diffeological spaces, $Y^X$ has for its underlying set
the set of smooth maps $X\to Y$; and a map $U\to Y^X$ is declared to
be a smooth plot if its transpose $U\times X \to Y$ is smooth.  The
inclusion into the smooth topos preserves the cartesian closed
structure.

     For any smooth manifold $M$, we have in particular a diffeology on
     $\ci (M) = {\bf R}^M$, namely a map $g:U\to \ci (M)$ is declared to be
     a smooth plot iff its transpose $U\times M \to {\bf R}$ is smooth.

     \medskip
     
 Topological vector spaces $V$ carry a canonical diffeology: a plot 
 $f: U\to V$ is declared to be smooth if for every continuous linear 
 functional $\phi :V \to {\bf R}$, $\phi \circ f : U \to {\bf R}$ is 
 smooth in the standard sense of multivariable calculus.
 
 A {\em convenient vector space}  is a toplogical vector 
 space  with certain properties,
cf.\ \cite{FK} 2.6.3; one of these is: if $\phi :V \to {\bf R}$ is a linear
functional, which is smooth with respect to the diffeological
structures on $V$ and ${\bf R}$, then $\phi $ is continuous.  (So the 
set of continous linear functionals is
closed under a certain obvious Galois correspondence between linear
functionals $V\to {\bf R}$, on the one side, and plots ${\bf R}^n
\supseteq U\to V$ on the other.) -- We use ``CVS'' as a shorthand for the
phrase ``convenient vector space''.

\medskip

{\bf Note:} Besides the category  of convenient
vector spaces as a full subcategory of the category of topological 
vector spaces, we shall have occasion
to consider another much larger cate\-gory $\underline{Con }^{\infty}$ of
convenient vector spaces; it has the same objects, but with {\em all}
smooth maps in between them, not just the smooth {\em linear} ones.
The category $\underline{Con}^{\infty}$ is a full subcategory of the category
$\underline{Diff}$ of diffeological spaces. 

\medskip

For convenient vector 
spaces, a map $f:X\to Y$ is plot smooth iff it is {\em scalarwise} smooth, 
meaning that $\phi \circ f : X\to {\bf R}$ is smooth for any $\phi 
\in Y'$ (where $Y'$ is the set of continuous (=smooth) linear 
functionals $Y\to {\bf R}$).

\medskip

Convenient vector spaces $X$ have the following completeness property:
given a smooth curve $g: U\to X$ (i.e.\ a smooth plot, where $U$ an open 
interval in  ${\bf R}$),
there is a unique function $g' :U\to X$ which is derivative of $g$ in
the scalarwise sense that $(\phi \circ g)' = \phi \circ g'$ for all $\phi \in
X'$; and this $g'$ is itself smooth.  -- More generally, if $U
\subseteq {\bf R}^n$ is open, and $g: U\to X$ is a smooth plot, then
partial derivatives $g^{\alpha}$ of $g$ exist, in the scalarwise sense; and
they are smooth (in particular, they are continuous).  Here $\alpha$ is a
multi-index; and to say that $g^{\alpha}$ is an iterated partial
derivative of $g$, in the {\em scalarwise} sense, is to say: for each $\phi
\in X'$, $\phi \circ g$ has an $\alpha$th iterated derivative, and
$(\phi \circ g)^{\alpha}= \phi \circ g^{\alpha}$.

Conversely, if $g:U \to X$ has the property that scalarwise iterated
partial derivatives $g^{\alpha}$ exist, and are scalarwise continuous,
then for each $\phi \in X'$, $\phi \circ g$ has iterated partial
derivatives, and they are continuous, since the $g^{\alpha}$ were
assumed to be so; so $\phi \circ g$ is smooth, and therefore $g$
itself is smooth.

\medskip

     For $i:X\to Y$ be a smooth linear map between convenient vector spaces.
Then $i$ preserves differentiation of smooth plots $U\to X$, in an
obvious sense.
For instance, if $f: U \to X$ is a smooth curve, i.e.\ $U\subseteq
{\bf R}$ an open interval, then for any $t_0 \in U$,
$$ (i\circ f)'(t_0 ) = i (f' (t_0 )).$$
For, it suffices to test this with the elements $\psi \in Y'$. If
$\psi \in Y'$, then $\psi \circ i \in X'$ since $i$ is smooth and
linear, and the result then follows by definition of being a scalarwise
derivative in $X$.

\medskip

     We don't know at present whether generally scalarwise smooth curves ${\bf
     R}_{\geq 0}\to X$ similarly have ``scalarwise derivatives'' in the
     endpoint $0$ (unless $X$ is ${\bf R}$, say, where the result follows
     from Theorem \ref{p33}).  This prompts us to make a definition.

\medskip

\begin{defn}  Call a map $g: {\bf R}_{\geq 0} \to X$ {\em
strongly smooth} if for each natural number $n$, there exists a map
$g^{(n)}: {\bf R}_{\geq 0} \to X$ such that for each $\phi \in X'$,
$\phi \circ g$ is $n$ times differentiable with $(\phi \circ g)^{(n)} =
\phi \circ g^{(n)}$.  \label{strong}\end{defn}

If $g: {\bf R}_{\geq 0} \to X$ is strongly smooth (with $X$ a CVS), it
is scalarwise smooth, in the sense that $\phi \circ g :
{\bf R}_{\geq 0} \to {\bf R}$ satisfies the condition (4) of Theorem
\ref{p33}, for every $\phi \in X'$. Hence, by the Theorem, $\phi \circ g$
is also  smooth for every $\phi \in X'$, and by definition of the
diffeology on the CVS $X$, $g$ itself is smooth.

So ``strongly smooth implies  smooth'', for maps ${\bf R}_{\geq 0} \to X$.

     \medskip

     (Another aspect of the completeness of convenient vector spaces is: if
      $U$ is an open interval, and $u_0 \in U$, there is a unique smooth
      primitive $G$ ($G' =g$) of $g$, with $G(u_0 )=0$.  This is the basis
      for constructing ``Hadamard remainders'' with values in a CVS, and
      hence for the comparisons of the present Section 5.)

      \medskip

Pointers to the literature: Convenient Vector Spaces were
introduced by Fr\"{o}licher and Kriegl, an exposition is in
\cite{FK}; diffeological spaces (cf.\  \cite{Souriau}) seem to have been
invented and re-invented with small variations and with different
names several times, because they seem not to be really admitted into
mainstream functional analysis.  One early reference is Chen's
\cite{Chen}, were a variation of the theme, under the name
``Differential Space'' is introduced.  Convenient Vector Spaces are
put into the context of diffeological spaces in \cite{Nel}.  A recent
account is given in \cite{Torre}, where also a comparison with
convenient vector spaces is presented.

The category $\underline{C}^{\infty}$ of {\em smooth spaces}, \cite{FK} 1.4.1,
is a full subcategory of the category of diffeological spaces, but
it does not enter directly in our exposition.

\section{The basic vector spaces of distribution theory; test plots}

Let $M$ be a smooth (paracompact) manifold $M$ (we shall here be interested
in ${\bf R}^m$, only).  Distribution theory starts out with the vector
space $C^{\infty}(M)$ of smooth real valued functions on $M$, and the
linear subspace ${\cal D}(M)\subseteq C^{\infty}(M)$ consisting of
functions with compact support (${\cal D}(M)$ is the ``space of {\em test
functions}).  The topology relevant for distribution theory is
described (in terms of convergence of sequences) in \cite{Schw1}, p.
79 and 108, respectively.  Note that the topology on ${\cal D}(M)$ is finer
than
the one induced from the topology on $C^{\infty}(M)$.  The sheaf
semantics which we shall consider in Section 7 will justify the choice
of these topologies.

We shall describe the diffeological structure, arising from the
topology on ${\cal D}(M)$, and utilize the fact (\cite{FK}, Remark 3.5)
that it is a  convenient
vector space.

      We cover $M$ by an increasing sequence
$K_b$ of compact subsets, $M=\cup K_b$; the notions that we now describe
are independent of the choice of these $K_b$.  For $M={\bf R}^n$, we would
typically take $K_b = \{x\in {\bf R}^n \mid |x| \leq b \}$, $b\in {\bf N}$.

Consider a smooth map $f: U\times {M} \to {\bf R}$, where $U$ is an
open subset of some ${\bf R}^n$.  We say that it is of {\em uniformly bounded
support} if there exists $b$ so that $$f(u,x)=0 \mbox{ for all } u\in
U \mbox{ and all } x \mbox{
with } x\notin K_b$$
We say that $f$ is  {\em locally} of uniformly bounded support
(``l.u.b.s.'') if $U$
can be covered by open subsets $U_i$ such that for each $i$, the
restriction of $f$ to $U_i \times M$ is of uniformly bounded
support.  (We may  use the phrase ``$f$ is
l.u.b.s., {\em locally in the variable} $u\in U$'') - Equivalently, we say
$f$ is
of uniform bounded support {\em at} $u\in U$ if there is an open
neighbourhood $U'$ around $u$ such that the restriction of $f$ to $U'
\times M$ is of uniformly bounded support; and $f$ is l.u.b.s.\ if it
for each $u$ is of uniformly bounded support at $u$.  (For yet another
description of the notion, see Lemma \ref{lubs2} below.)

\medskip

     We let $\hat{f}$ denote the transpose of $f$, so $\hat{f}: U\to
C^{\infty}(M)$.

\begin{thm} Let $f:U\times M\to {\bf R}$ be smooth, and pointwise of
bounded support (so that $\hat{f}$ factors through ${\cal D}(M)$).
Then t.f.a.e.:

1) $f$ is locally of uniformly bounded support

2) $\hat{f}: U \to {\cal D}(M)$ is continuous.
\label{one}\end{thm}
We may use the term {\em test plot} for functions $f$ satisfying the
conditions of the Proposition.  Pointwise, they are test {\em functions}
in the sense of distribution theory.

\medskip

\noindent {\bf Proof} of the Theorem.  We first prove that 1) implies 2).
Since the question is local in $U$, we may assume that $f$ is of
uniformly bounded support, i.e.\ there exists a compact $K\subseteq
M$ so that $f(t,x)=0$ for $x\notin K$ and all $t$.  The same $K$ applies then
to all the iterated partial derivatives $f_{\alpha}$ of $f$ in the
$M$-directions ($\alpha$ denoting some multi-index).  So $f$ and all
the $f_{\alpha}$ factor through ${\cal D}_K$, the subset of $\ci (M)$
of functions vanishing outside $K$.  Now to say that $\hat{f}:U \to {\cal
D}_K$ is
continuous is by definition of the topology on ${\cal D}_K$ equivalent
to saying that for each $\alpha$, $(f_{\alpha})^{\hat{}}$ is continuous
as a map into ${\bf R}^K$, the space of continuous maps $K\to {\bf
R}$, with the topology of uniform convergence.  This topology is the
categorical exponent ( = compact open topology) (cf.  \cite{Kelley}
Ch.\ 7 Thm.\ 11), which implies that $(f_{\alpha})^{\hat{}}: U\to {\bf
R}^K$ is continuous iff $f_{\alpha} :U\times K\to {\bf R}$ is
continuous, iff $f_{\alpha} :U\times M\to {\bf R}$ is continuous.  But
$f_{\alpha}$ is indeed continuous, by the smoothness assumption on
$f$.  So $\hat{f}:U \to {\cal D}(M)$ is continuous.

        For proving that 2) implies 1), we prove that if not 1), then not 2),
i.e.\
        we consider a function $f : U\times M \to {\bf R}$ which is smooth and
        of pointwise bounded support, but not l.u.b.s.  Then there is a $t_0
        \in U$ and a sequence $t_k \to t_0$, as well as  a sequence $x_k \in M
        \setminus K_k$ with $f(t_k , x_k) \neq 0$, denote this number $c_k$.
        Let $N$ be a number so that the support of $f(t_0, -)$ is contained in
        $K_N$.  We consider the (non-linear) functional $T: {\cal D}(M) \to
        {\bf R}$ given by
$$g\mapsto \sum _{n=N} ^{\infty} c_n ^{-2} g(x_n )^2.$$
Note that for $g$ of compact support, this sum is finite, since the
$x_n$'s ``tend to infinity''.  Also, the functional ${\cal D}(M) \to
{\bf R}$ is continuous; for the topology on ${\cal D}(M)$ is the
inductive limit of the topology ${\cal D}(K_k )$, and the restriction
of $T$ to this subspace equals a {\em finite} algebraic combination of
the Dirac distributions.  Now it is easy to see that $T$ takes $f(t_0,
-)$ to $0$, by the choice of $N$, whereas
$T$ applied to $f(t_k , -)$ for $k>N$ yields a sum of non-negative
terms, one of which has value 1, namely the one with index $k$, which is
$c_k ^{-2} f(t_k , x_k )^2 =1$.  So $T\circ \hat{f}$ is not
continuous, hence $\hat{f}$ is not continuous.

This proves the Theorem.

\medskip

It has the following Corollary:

\begin{thm} Let $f: U \times M \to {\bf R}$ be smooth and of pointwise bounded
support ($U$ an open subset of some ${\bf R}^n$).  Then t.f.a.e.:

1) $f$ is  locally of uniformly bounded support

2) $\hat{f}:
U \to {\cal D}(M)$ is smooth.

\label{lubs}\end{thm}
Recall that assertion 2) means ``in the scalarwise sense'', i.e.\ 
$\phi \circ \hat{f}$ is smooth for any continuous linear functional, 
i.e.\ for any distribution  
$\phi$.

{\bf Proof.} The implication 2) implies 1) is a consequence of
Theorem \ref{one}, since smoothness implies continuity.  Conversely, assume
1), i.e.\ assume $f$ is smooth and l.u.b.s.  Then we also have that
$\partial ^{\alpha} f / \partial t^{\alpha}$ is smooth (iterated
partial derivative in the $U$-directions, $\alpha$ a multi-index) and
l.u.b.s., and so its transpose is a continuous maps $U\to {\cal
D}(M)$, by Theorem \ref{one}; it serves as scalarwise iterated partial
derivative.  (This is an entirely classical statement; we could not
find an explicit reference, so we sketch a proof.  The continuous
linear functionals that define what ``scalarwise'' means are by
definition the distributions on $M$.  For the case where $M$ and $U$
both are ${\bf R}$, it is thus the assertion that for any distribution
$T$, if $f(t,s): {\bf R}\times {\bf R} \to {\bf R}$ is smooth, and,
locally in the variable $t$, of uniformly bounded support, then
$t\mapsto T([s\mapsto f(t,s)])$ has a $t$-derivative which is given by
$t\mapsto T([s\mapsto \partial f (t,s)/\partial t])$, i.e.\ ``one can
differentiate under the distribution sign''.  Now since the desired
conclusion is of local nature in $t$, we may w.l.o.g.  assume that $f$
is of uniformly bounded support, and then we may modify $T$ so as to
have compact support also.  Then $T$ may be represented by a finite
sum of ``derivatives of continuous functions'' (cf.  \cite{Schw2}
Thm.\ 26), so the assertion of ``differentiating under the
distribution sign'' becomes essentially the assertion that you may
differentiate under the integration sign, which is possible due to the
compactness of the support.) \medskip

The standard vector space of distributions ${\cal D}'(M)$ is, in
diffeological terms, the linear subspace of the diffeological space
${\bf R}^{{\cal D}(M)}$ consisting of the {\em linear} smooth maps
${\cal D}(M) \to {\bf R}$.  A map $U\to {\cal D}'(M)$ is smooth iff it
is smooth as a map into ${\bf R}^{{\cal D}(M)}$; this defines a
diffeology on ${\cal D}'(M)$.  With this diffeology, ${\cal D}'(M)$,
too, is convenient.

The diffeology/convenient vector space structure on ${\cal D}(M)$
corresponds to its standard locally convex topology, so that a linear
functional ${\cal D}(M) \to {\bf R}$ is diffeological iff it is
continuous.  So the vector space of distributions ${\cal D}'(M)$ (as
an abstract vector space) is the same in both contexts.

We have
\begin{thm} The convenient vector space ${\cal D}(M)$ is reflexive in
the CVS sense: the canonical ${\cal D}(M) \to {\cal D}''(M)$ is an
isomorphism.
\label{refl}\end{thm} We presume that this result is well known among experts,
but it is not explicitly stated in \cite{FK}, say.  It is known that,
as a locally convex topological vector space, ${\cal D}(M)$ is
reflexive, with respect to the so called strong topology on dual
spaces, cf.  \cite{Schw2} Theorem XIV.  Now \cite{FK} has a general
Theorem, comparing reflexivity in various categories of vector spaces
(locally convex, convenient, bornological, \ldots), namely Theorem
5.4.6.  By this Theorem, convenient reflexivity follows from strong
reflexivity, provide that the strong dual (in our case ${\cal D}'(M)$
with its strong topology) is furthermore bornological.  And this is
known to be so, cf.  e.g.\ \cite{Horvath} Example 3.16.2.

\section{Functions as distributions}  Any sufficiently nice function
$f:{\bf R}^n \to {\bf R}$ gives
rise to a distribution $i(f) \in {\cal D}' ({\bf R}^n )$ in the standard way
``by integration over
${\bf R}^n$''
$$<i(f), \phi > := \int _{{\bf R}^n} f(s)\cdot \phi (s) \; ds.$$
This also applies if ${\bf R}^n$ is replaced by another smooth
manifold $M$ equipped with a suitable measure. For simplicity of
notation, we write $M$ for ${\bf R}^n$ in the following.
-- All smooth functions $f:M \to {\bf R}$ are ``sufficiently nice'';
so we get a map (obviously linear)
\begin{equation}
i: \ci (M) \to {\cal D}' (M).
\label{i-map}\end{equation}
It is also easy to see that this map is injective.
\begin{thm} The map $i$ is smooth.
\label{ismooth}\end{thm} {\bf Proof.} Let $g: V\to \ci (M)$ be smooth, ($V$
an open subset of some ${\bf R}^n$), we have to see that $i\circ g : V
\to {\cal D}' (M)$ is smooth, which in turn means that its transpose
$$(i\circ g)\hat{} :V\times {\cal D}(M) \to {\bf R}$$
is smooth. So
consider a smooth plot $U \to V\times {\cal D}(M)$, given by a pair of
smooth maps $h: U\to V$ and $\hat{\Phi} :U\to {\cal D}(M)$. Here $U$
is again an open subset of some ${\bf R}^k$. Let us write $\hat{F}$ for
$g\circ h : U \to \ci (M)$. It is transpose of a map $F: U\times M
\to {\bf R}$. Also,  let us write $\Phi$ for the transpose of
$\hat{\Phi }$; thus $\Phi $ is a map
$$\Phi : U\times M \to {\bf R}$$
which is locally (in $U$) of uniformly bounded support, by Theorem \ref{lubs}.
We have to see that $(i\circ g)\hat{} \; \circ <h,\Phi >$ is smooth (in
the usual sense).  By unravelling the transpositions, one can easily
check that
$$(i\circ g)\hat{} \; \circ <h,\Phi > (t)= <i(F(t,-), \Phi (t,-)>$$
    The conclusion of the Theorem is thus the assertion that the composite
    map $U\to {\bf R}$ given by \begin{equation}t\mapsto \int _M F
    (t,s)\cdot \Phi (t,s) \; ds
\label{smooth?}\end{equation}
is smooth (in the standard sense of finite dimensional calculus). To
prove smoothness at $t_0 \in U$, we may find a neigbourhood $U'$ of
$t_0$ and a $b$ such that
$$\Phi (t,s) =0 \mbox{ if } t\in U' \mbox{ and } s\notin K_b,$$ because
$\Phi $ is l.u.b.s.  We thus have, for any $t\in U'$, that the
expression in (\ref{smooth?}) is $\int _{K_b}F (t,s)\cdot \Phi (t,s)
\; ds $, but since $K_b$ is compact, differentiation and other limits
in the variable $t$ may be taken inside the integration sign.

\medskip

Since $i:\ci (M) \to {\cal D}'(M)$ is smooth and linear, it preserves
differentiation. In particular, if $f: U \to \ci (M)$ is a smooth
curve, and $t_0 \in U$, we have that
$(i\circ f)'(t_0 ) = i (f' (t_0 ))$. However, $f'$ is
explicitly calculated in terms of the partial derivative of the
transpose $\hat{f} : U\times M \to {\bf R}$, namely as the function
$s\mapsto \partial f(t,s) /\partial t \mid _{(t_0 ,s)}$.
This is the reason that ordinary \mbox{(evolution-)} diffe\-rential equations
for curves $f:U\to {\cal D}'(M)$ manifest themselves as {\em partial}
differential equations, as soon as the values of $f$ are distributions
represented by smooth functions.

\section{Smoothness of heat kernel}

We consider the heat equation on the line,
$$\partial f/\partial t (t,x)= \partial ^2 f /\partial x^2 .$$

Recall that the classical distribution solution of this equation, having
$\delta (0)$ as initial distribution, is   the
map
$$K: {\bf R}_{\geq 0}\lra {\cal D}'({\bf R})$$
whose value at $t\geq 0$ is the distribution $<K(t),\;\_\;>$ given on
a test function $\phi$ by
\begin{equation} <K(t),\phi> = \left\{ \begin{array}{ll} \int^{\infty}_{\-
\infty}
e^{-s^2/4t}/\sqrt{4\pi t}\; \phi (s)\; ds & \mbox{if $t>0$} \\
\phi(0) & \mbox{if $t=0$}
\end{array}
\right.
\label{gauss}\end{equation}

The present Section is devoted to proving the strong
smoothness (Definition \ref{strong}), and hence also the diffeological
smoothness, of $K$.

\begin{thm}The function $K: {\bf R}_{\geq 0} \to {\cal D}' ({\bf R})$
is strongly smooth.
\label{strongsmooth}\end{thm}
{\bf Proof.} We have to produce for each $n$  a map $K^{(n)}:{\bf
R}_{\geq 0} \to {\cal D}' ({\bf R})$ which will serve as the $n$'th
scalarwise derivative of $K$.  The map $\Delta ^n \circ K$ will do.
For, consider a smooth linear $\rho : {\cal D}'({\bf R})\to {\bf R}$.
So $\rho \in {\cal D}'' ({\bf R})$, but by CVS-reflexivity of ${\cal
D}({\bf R})$ (Theorem \ref{refl}), $\rho $ is of the form $\rho (T) =
<T, \phi >$ for a unique test function $\phi \in {\cal D}({\bf R})$.

So it suffices to prove, for each test function $\phi$,  that $<(\Delta
^n \circ K)(t), \phi >$,
as a function of $t\in {\bf R}_{\geq 0}$, is the $n$'th derivative of
$<K(t), \phi >$. Now, $$<(\Delta
^n \circ K)(t), \phi > =<K(t), \phi ^{(2n)}>$$ (recalling
$\Delta = (-)''$, and the differentiation of distributions); so this
reduces the problem to proving the following:

\begin{prop}\label{hs2}
Let $\phi$ be a smooth function of compact support and let $\Phi:{\bf R}_{\geq
0}\lra {\bf R}$ be the function defined by
$$ \Phi(t)=<K(t),\phi>.$$
Then the function $\Phi$ is smooth and, furthermore, for all $t\geq 0$,
${\Phi}^{(n)}(t)={\phi}^{(2n)}(t)$.
\end{prop}
For $t>0$, this is well known:  for all $t>0,$ the $n$'th derivative of
$\Phi$ exists at
$t$, and
\begin{equation}{\Phi}^{(n)}(t)=<K(t), {\phi}^{(2n)}>,\label{first}
\end{equation} see p.\ 330 in \cite{Schw1}).  Next, we prove
\begin{lemma}
For every test function $\phi :{\bf R}\to {\bf R}$
$$lim_{t\ra 0^+}(1/t)[<K(t), \phi>-\; {\phi}(0)]=\phi''(0) $$
\label{second}\end{lemma}
{\bf Proof:} We first notice that, by Hadamard's Lemma, $\phi (x)=\phi
(0) + x\psi (x)$, for a unique smooth function $\psi$. (The function $\psi$
goes to $0$ when $x\ra +\infty$ or $x\ra -{\infty}$ since
$\psi(x)=(1/x)[\phi(x)-\phi(0)]$, but does not necessarily have
compact support.  Similarly for $\psi '$, $\psi ''$, etc.; but this
boundedness is enough to make the improper integrals convergent.) Note
that $\phi ''(0) = 2\psi ' (0)$.

We claim that for $t>0$ \begin{equation}<K(t), \phi >- \phi(0) = 2t<K(t),
\psi'>.\label{twostar}\end{equation} In fact, start from the right
hand side; we get the limit as $N_1 $ and $N_2 \to \infty$ of
$$2t\int_{-N_1}^{N_2 }1/\sqrt{4\pi t}\;
e^{-{x^2}/4t}\psi'(x)dx$$
which we integrate by parts to get
     $$(2t\cdot 1/\sqrt{4\pi t}\; e^{-{x^2}/4t}
\cdot \psi (x))|_{-N_1}^{N_2}+ \int_{-N_1}^{N_2}(1/\sqrt{4\pi
t})e^{-x^2/4t}x\psi (x)\;  dx .$$
Since $\psi$ is bounded, the first term here tends to 0 as $N_1, N_2
$ tend to $\infty$, so passing to the limit, and using $x\psi (x)= \phi
(x)-\phi (0)$, we get
$$2t <K(t) , \psi '> = \int _{-\infty}^{\infty} \frac{1}{\sqrt{4\pi
t}}e^{-x^2 /4t} x\psi (x) \; dx$$
$$=\int _{-\infty}^{\infty} \frac{1}{\sqrt{4\pi
t}}e^{-x^2 /4t} (\phi (x)-\phi (0) \; dx$$
$$=\int _{-\infty}^{\infty} \frac{1}{\sqrt{4\pi
t}}e^{-x^2 /4t} \phi (x)\; dx -\phi (0)$$
$$= <K(t), \phi > - <K(0), \phi>.$$
Now divide by $t$ and let $t\to 0$, recalling that $<K(t), \psi '> \to
\psi '(0)$ as $t\to 0$.

\medskip
\noindent {\bf Proof} of the Proposition: We first show that
\begin{equation}{\Phi}^{(n)}(0)={\phi}^{(2n)}(0)\label{induct}\end{equation}
by induction on $n$.
For $n=0$ this is by definition of $\Phi$. Assume that it is true for
$n$.
Then
$$\begin{array}{lll}
{\Phi}^{(n+1)}(0) & = & lim_{t\ra
0^+}(1/t)[{\Phi}^{(n)}(t)-{\Phi}^{(n)}(0)] \\
                  & = & lim_{t\ra 0^+}(1/t)[<K(t), \phi^{(2n)}>-\;
{\phi}^{(2n)}(0)] \\
                  & = & ({\phi}^{(2n)})''(0) \\
                  & = &  {\phi}^{(2(n+1))}(0)
\end{array}$$
In the passage from the first line to the second we have used the
induction hypothesis and  (\ref{first}), whereas to go to the third line
from the second we have used Lemma \ref{second} with $\phi^{(2n)}.$

\medskip

We noted already that for a map with domain ${\bf R}$ or ${\bf
R}_{\geq 0}$, and with values in a CVS, strong smoothness implies
scalarwise (equivalently, diffeological) smoothness. So Theorem
\ref{strongsmooth} has as Corollary:

\begin{thm}The function $K: {\bf R}_{\geq 0} \to {\cal D}' ({\bf R})$ is
smooth (in the diffeological sense).
\label{fin}\end{thm}

   This Theorem, in turn, translates by passing to the transpose map, the
   following result, formulated in entirely classical and elementary-calculus
   terms.  (In one of the preliminary versions of the present
   paper, we gave an elementary proof of it, and hence of Theorem
   \ref{fin}, without resorting to reflexivity of ${\cal D}(M) $.)
   \begin{thm} Let $\phi : U \times {\bf R} \to {\bf R}$ be a test plot
   (i.e.\ smooth, and locally (on $U$) of uniformly bounded support).
   Then the function $\Phi :U \times {\bf R}_{\geq 0} \to {\bf R}$
     defined by $$\Phi (u,t) :=<K(t),\phi (u,-)>$$
is smooth .
\label{param}\end{thm}

\section{Ideals and differential operators}

Let $x\in {\bf R}^n$. By a {\em differential operator supported at} $x$, we
understand a map $d: \ci ({\bf R}^n )\to {\bf R}$ which is a linear
combination of operators $f\mapsto \partial ^{|\alpha |} f / \partial
t^{\alpha } (x)$, where $\alpha$ is a multi-index and $t = (t_1 ,
\ldots ,t_n )$.  (The notion can be defined in a coordinate free way;
it is actually the same as a distribution with point-support.) In
particular, $d$ is linear.

Any such $d$ defines, because of its explicit form, for each CVS $Y$ a linear
$d_Y :\ci ({\bf R}^n , Y) \to Y$ with the property that for $f:{\bf
R}^n \to Y$
$$ d(\phi \circ f) = \phi (d_Y (f))$$
for all $\phi \in Y'$.  The maps $d_Y$ are natural in $Y$ w.r.to
smooth {\em linear} maps:

\begin{prop} If $F: Y\to X$ is a smooth linear map, then for any
differential operator $d$, and any $f\in \ci ({\bf R}^n , Y)$, $d_X
(F\circ f) = F(d_Y (f))$
\label{i-preserves}\end{prop}
{\bf Proof.} It suffices to test with an arbitrary $\phi \in X'$; by
replacing $F$ by $\phi \circ F$, this reduces the problem to the case
where the dodomain $X$ is ${\bf R}$, and here, the result follows from
the very characterization of $Y$-valued derivatives in ``scalarwise'' terms.

\medskip
Let us also note that ``partial derivatives are transposable''.  For
simplicity, we state it for functions in two variables $s,t$ only:

\begin{prop} Let $f(s,t): {\bf R}^2 \to Y$ be a smooth function with
values in a CVS $Y$.  Then the $\partial f (s , t) /\partial s $ is
smooth in $s,t$, and its transpose is the derivative $(\hat{f})'(s)$ of the
transposed function $\hat{f}: {\bf R} \to \ci ({\bf R}, Y)$.
\label{transposability} \end{prop}
{\bf Proof.} The function $(\hat{f})'(s)$ exists and is smooth, and
characterized  in terms of the smooth linear
     functionals on $\ci ({\bf R}, Y)$.  But among these are those of the
     form (for $t\in {\bf R}$) \begin{diagram} \ci ({\bf R}, Y)&\rTo ^{ev_t} &
Y &\rTo
     ^{\psi} &{\bf R},\end{diagram} and these are enough to recognize the
     transpose of $\hat{f}'(s)$ as $\partial f (s , t) $, for each $t$.

\medskip

Let $I \subseteq \Ci $ be an ideal.  For each CVS $Y$, (in fact for
any dualized vector space $(Y,Y')$) we define two linear subspaces of $\ci
({\bf R}^l , Y )$, the ``weak'' and the ``strong'' $I(Y)$, denoted
$I_w (Y)$ and $I_s (Y)$, respectively.  To say that $f: N \to Y$ is in
$I_w (Y)$ is to say that for every $\phi \in Y'$, $\phi \circ f \in
I$; and to say that $f: N \to Y$ is in $I_s (Y)$ is to say that $f$
may be written
$$f(s)= \sum h_i (s) k_i (s),$$
with the $h_i$'s scalar valued functions belonging to $I$, and the
$k_i$'s smooth $Y$-valued functions. It is clear that
$I_s (Y)\subseteq I_w (Y)$.  We are interested in when the converse
implication holds.

A main result in \cite{CSF} (Theorem 2.11) says that this is the case
for the ideal ${\cal M}^r \subseteq C^{\infty} ({\bf R}^l)$ of
functions vanishing to order $r$ at $0$.  In \cite{CA} (Proposition
1), we generalized this to any proper ideal $I \subseteq C^{\infty} ({\bf
R}^l)$ which contains an ideal ${\cal M}^r$.  We call such ideals {\em
Weil ideals}; they are of finite codimension, and the algebra $
C^{\infty} ({\bf R}^l)/I$ is a Weil algebra (in the sense of
\cite{SDG} or \cite{MR}, say); and any Weil algebra arises this way.
(Note that a Weil ideal is contained in ${\cal M}$, since the only
maximal ideal containing ${\cal M}^r$ is ${\cal M}$.  So if $f \in I$,
$f(0)=0$.)

\medskip

We shall generalize this result further to what we call ``semi-Weil
ideals'' $J$, and at the same time provide a simpler proof of the
result quoted  from \cite{CA}.

\medskip

If $I \subseteq \cz (N)$ is an ideal and if $p: P\to N$ is a smooth map
($P$ and $N$ manifolds), we get an ideal $p^* (I) \subseteq \cz (P)$
consisting of functions $f:P\to {\bf R}$ which can be written $\sum
(h_i \circ p) \cdot k_i$ with the $h_i$'s in $I$ (and the $k_i$'s in
$\ci (P)$).  This is clearly a ``transitive'' construction, in an evident
sense, $q^* (p^* (I))= (p\circ q)^* (I)$.  On the other hand, since
$\cz (M)$ is a convenient vector space, we may consider $I_s (\cz
(M))\subseteq \cz (N , \cz (M))$.  Under the isomorphism $\cz (N , \cz
(M)) \cong \cz (N \times M)$ it is clear that $I_s (\cz (M))$
corresponds to $p^* (I)$, where $p: N\times M \to N$ denotes the
projection.

If $I$ is a Weil ideal $\subseteq \Cil $, and
$p: {\bf R}^{l+k}\to {\bf R}^l$ the projection, we get by the above
procedure an ideal $J=p^* (I)$ in $\cz ({\bf R}^{l+k})$, and ideals $J$
of this form, we call {\em semi-Weil ideals}.

The basis of monomials $s^{\alpha}$ (where $\alpha$ is a multi-index of
order $<r$) for $C^{\infty}({\bf
R}^l)/({\cal M}^r)$ gives rise to a dual basis for the
linear dual $(C^{\infty}({\bf
R}^l)/{\cal M}^r))^*$, and this dual basis consists of differential
operators supported at $0$, $$f \mapsto \frac{\partial ^{\alpha}f
(0)}{|\alpha |!  \partial s^{\alpha}}.$$   So $f\in {\cal M}^r$ iff
$\frac{\partial ^{\alpha}f (0)}{
s^{\alpha}}=0$ for such multi-indices $\alpha$.

\medskip

We consider functions $f(s,t): {\bf R}^{l+k} \to Y$, where $Y$ is a
CVS; $s$ denotes a variable ranging over ${\bf R}^l$ and $t$ a
variable ranging over ${\bf R}^k$. We then have
\begin{prop}Let $f:{\bf R}^{l+k} \to Y$ be a smooth function. Then
$$f\in (p^* ({\cal M}^r))_w (Y)$$ if and only if $$\frac{\partial
^{\alpha}f (0,t)}{
\partial s^{\alpha}}=0 \mbox{ for all } \alpha \mbox{ with } |\alpha | <r
\mbox{ and
all }t.$$
\label{multi}\end{prop}
{\bf Proof.} The $Y$-valued partial derivatives here are determined
scalarwise, i.e.  determined by testing with the $\phi \in Y'$, and since
these $\phi$ are linear, the problem immediately reduces to the case
of $Y={\bf R}$, i.e.  to the assertion $f(s,t) \in p^* ({\cal M}^r )$
iff $\frac{\partial ^{\alpha}f (0,t)}{ s^{\alpha}}=0 $ for all
$\alpha$ with $|\alpha | <r$ and all $t$.  This is well known (or can
be deduced from Theorem 2.11 in \cite{CSF}, by passing to the
transpose function $\hat{f}:{\bf R}^l \to \ci ({\bf R}^k )$).

\medskip

The following is now a Corollary of Theorem 2.11 in \cite{CSF}:
\begin{prop}For any CVS $Y$, we have $(p^* ({\cal M}^r))_w (Y) =
(p^* ({\cal M}^r ))_s (Y).$
\label{fu}\end{prop}
{\bf Proof.} It suffices to prove the inclusion $\subseteq$. If $f$
is in the left hand side, it satisfies the equational conditions of
Proposition \ref{multi}, but then its transpose $\hat{f}:{\bf R}^l
\to \ci ({\bf R}^k , Y)$ has  $\frac{\partial
^{\alpha} \hat{f} (0)}{
s^{\alpha}}=0 $ for all $\alpha$ with $|\alpha | <r$. Now we apply
Proposition \ref{multi} again, this time for the CVS $\ci ({\bf R}^k ,
Y)$, and with no $p^*$ involved, and conclude $\hat{f} \in ({\cal
M}^r)_w (\ci ({\bf R}^k ,
Y) )$. Then, by the Theorem quoted, $\hat{f} \in ({\cal
M}^r)_s (\ci ({\bf R}^k ,
Y) )$ (strong instead of weak), and this in turn implies that $f \in
(p^* ({\cal M}^r))_s (Y)$, proving the Proposition.

\medskip

Consider a Weil ideal $I$ i.e.\ an ideal $I\subseteq \Ci $
containing some ${\cal M}^r$.
There is a (finite) basis $A$ for the dual vector space $(\Ci /{\cal M}^r )^*$
consisting of differential operators $D^{\alpha}$ at $0$ (with ${\cal M}^r$
the common nullspace of these).  Since $(\Ci /I)^* \subseteq (\Ci
/{\cal M}^r )^*$, we may, by suitable change of basis, organize
ourselves so that the basis $A$ for $(\Ci /{\cal M}^r)^*$ contains a
subset $B$ which is a basis for $(\Ci /I)^*$.  It follows that $I$ is
the common null space of the collection $B$ of differential operators.

The dual basis $\hat{A}$ for $\Ci /{\cal M}^r$ consists (modulo ${\cal
M}^r$) of polynomials $h_{\alpha}$ of degree $<r$, ($\alpha \in A$).
The fact that the bases $A$ and $A'$ are dual implies that for any $f\in \Ci$,
$$f(s) \equiv \sum _{\alpha} D^{\alpha}f \cdot h_{\alpha} (s),$$
mod ${\cal M}^r$ (as functions of $s \in {\bf R}^l$). If now $f\in
I$, the terms
$D^{\alpha}f$ vanish for $\alpha \in B$.  With $A-B$ as index set for
the index $\gamma$, we therefore have

\begin{prop}Given a Weil-ideal $I \in \Ci $ containing ${\cal M}^r$.
There is a finite family
of differential operators $D^{\gamma}$ and a family of polynomials
$h_{\gamma}(s)$ in $s\in {\bf R}^l$ so that for any $f \in I $,
$$f(s) -\sum _{\gamma} D^{\gamma}f \cdot h_{\gamma}(s) \in {\cal
M}^r .$$
\label{reduc}\end{prop}

(If for instance $I= {\cal M}^{r-1} \supseteq {\cal M}^r$, the
$h_{\gamma}$'s may be taken to be the monomials $s^{\alpha}$, where
$\alpha$ ranges over multi-indices with $|\alpha | =r$.)

\medskip

Because differentiation of functions ${\bf R}^l \to Y$ (with $Y$ a CVS)
makes sense, and because of the explicit way (in terms of
$D^{\gamma}$'s) in which functions in $I$ get transformed into
functions in ${\cal M}^r$, this Proposition immediately extends to
functions ${\bf R}^l \to Y$; let $I$ and $h_{\gamma}$ be as above, and
let the $D^{\gamma}$ denote the $Y$-valued differential operators
corresponding to the ${\bf R}$-valued $D^{\gamma}$'s considered.

\begin{prop} For any $f\in I_w (Y)$, the difference
$$f(s) -\sum _{\gamma} D^{\gamma}f \cdot h_{\gamma}(s) $$ belongs to
${\cal
M}^r _w (Y)$ (which equals $ {\cal
M}^r _s (Y))$ by the Theorem \cite{CSF} 2.11 quoted).
\label{y1}\end{prop}
{\bf Proof.} We test with arbitrary $\phi \in Y'$; since $\phi$ is
linear, and since $\phi$ commutes with differentiation, the result
follows by applying the result of the previous Proposition to the
smooth function $\phi \circ f $, which is in $I$ by assumption.

\medskip

Now let $J$ denote the semi-Weil ideal  $p^* I \subseteq C^{\infty} ({\bf
R}^{l+k})$ given by the Weil ideal $I \subseteq \Cil$.  Then

\begin{prop} Let $f:{\bf R}^{l+k} \to {\bf R}$ be a function in $J$. Then
$$f(s,t) - \sum (D^{\gamma}f)(t) \cdot h_{\gamma} (s)$$
is in $p^* ({\cal M}^r )$ (where $p:{\bf R}^{l+k} \to {\bf R}^l$ is
the projection).
\label{y2} \end{prop} (Here, $s$ and $t$ denote variables ranging over ${\bf
R}^l$ and ${\bf R}^k$, respectively.  The differential operators
$D^{\gamma}$ operate in the $s$-variable and then $s=0$ is
substituted, so a function $D^{\gamma}f$ of $t$ remains, as
indicated.)

\medskip

{\bf Proof.} We pass to the transpose function $\hat{f} : {\bf R}^l
\to Y$, where $Y$ is the CVS $\ci ({\bf R}^k )$.  To say $f\in J$ is
equivalent to saying $\hat{f} \in I_s (\ci ({\bf R}^k ))$, in
particular $\hat{f} \in I_w (\ci ({\bf R}^k ))$, and so Proposition
\ref{y1} may be applied, reducing $\hat{f}$ to ${\cal M}^r _s (\ci
({\bf R}^k )$, which by transposition corresponds to $p^* ({\cal M}^r
)$.  This proves the Proposition.

\medskip

We generalize this further to the case of functions with values in a CVS $Y$.

\begin{prop} Let $J \subseteq \ci ({\bf R}^{l+k})$ be the semi-Weil
ideal given by the Weil ideal $I \in \Ci $.  Let $g(s,t) \in J_w (Y)$. Then
\begin{equation}g(s,t)- \sum (D^{\gamma} g)(t)\cdot h_{\gamma }(s)
\label{cruc}\end{equation}
  is in $(p^* ({\cal M}^r ))_w (Y)$ (hence, by
Proposition \ref{fu}, in $(p^* ({\cal M}^r ))_s (Y)$).
\label{y4}\end{prop}
{\bf Proof.} Testing with $\phi \in Y'$ reduces the problem to showing that
\begin{equation}\phi (g(s,t)) - \sum D^{\gamma} (\phi \circ g ) (t) \cdot
h_{\gamma}(s)
\end{equation}
is in $p^* ({\cal M}^r )$, but this follows from
Proposition \ref{y2}, applied to $f= \phi \circ g$.

\begin{thm}If $J$ is a semi-Weil ideal, and $Y$ a CVS, $J_s (Y) = J_w
(Y)$ (as linear subspaces of $C^{\infty}({\bf R}^{l+k}, Y)$).
\label{swsw}\end{thm}
{\bf Proof.} Let $g=g(s,t)$,  $g:{\bf R}^l \times {\bf R}^k \to Y$,
be a map  in
$J_w (Y)$. Since the $h_{\gamma }(s)$ are in $I$, the sum
$\sum (D^{\gamma} g)(t)\cdot h_{\gamma }(s)$ in (\ref{cruc}) is in
$J_s (Y)$.  The whole expression in (\ref{cruc}) is in $(p^* ({\cal M}^r
))_w (Y)$, by Proposition \ref{y4}, and hence, by Proposition
\ref{fu}, in $(p^* ({\cal M}^r ))_s (Y)$ which in turn is contained in
$ J_s (Y)$.  This proves the Theorem.

\medskip

   From now on, we write $J(Y)$ instead of $J_w (Y)$ or $J_s (Y)$, in
case $J$ is a semi-Weil ideal and $Y$ a CVS; for, they agree, by the Theorem.

\medskip

We now discuss the description of semi-Weil ideals in terms of
differential operators.

    If $I\subseteq \ci ({\bf R}^n )$ is an ideal which is the null
space of a family of differential operators $\{d^{\beta} \mid \; \beta
\in B\}$ (not necessarily supported at the same $x\in {\bf R}^n$),
then it follows from Proposition \ref{i-preserves} that $I_w (Y)
\subseteq \ci ({\bf R}^n , Y)$ is the null space of the family of the
$d^{\beta}_Y$.

If $I$ is a Weil ideal in $\ci ({\bf R }^l )$, null
space of a finite family $\{ d^{\beta}\mid  \beta \in B\} $ of
differential operators supported at $0\in {\bf R}^l$, then $J
\subseteq \ci ({\bf R }^{l+k})$ is the null space of the (infinite)
family of differential operators $d^{\beta , x}$, $\beta \in B$, $x\in
{\bf R}^k$, where for a function $f(s,t) \in \ci ({\bf R }^{l+k})$,
$d^{\beta , x} (f)$ takes the relevant partial derivatives in the
$s$-directions, and then substitutes $0$ for $s$ and $x$ for $t$.

It follows that $J (Y)$, for $Y$ a CVS, may be described as the null
space of the $B\times {\bf R}^k$-indexed family of differential
operators $d^{\beta , x} _Y : \ci ({\bf R}^{l+k} , Y) \to Y$.

Also, it follows that under the transposition isomorphism
$\ci ({\bf R}^{l+k} ,Y)\cong \ci ({\bf R}^l , \ci ({\bf R}^k , Y))$,
the linear subspace $J(Y)$ on the left corresponds to the linear
subspace $I(\ci ({\bf R}^k , Y)$ on the right.

\medskip

    Let $I \subseteq {\bf R}^l$ be a Weil ideal, $I\supseteq
{\cal M}^r$.  Let $\{ D^{\beta}\mid  \beta \in B\} $ be a family of
differential
operators at $0$, of degree $<r$, forming a basis for $(\ci ({\bf R}^l )/I
)^*$.  Note that $B$ is a finite set.  Let the dual basis for $(\ci ({\bf
R}^l )/I$
     be represented by polynomials of degree $<r$, $\{ p_{\beta} (s) \mid
\beta \in
     B\}$.  Then we can construct a linear isomorphism
$$\ci ({\bf R}^l ,Y)/I(Y) \to \prod  _B Y,$$
by sending the class of $f:{\bf R}^l \to Y$ into the $B$-tuple
$D^{\beta}_Y (f)$.  Its inverse is given by sending a $B$-tuple
$y_{\beta} \in Y$ to $\sum _B p_{\beta} (s)\cdot y_{\beta}$.

It follows that for a semi-Weil ideal $J =p^* (I) \subseteq {\bf
R}^{l+k}$, as above,
\begin{equation}\ci ({\bf R}^{k+l} , Y)/J(Y) \cong \prod  _B
\ci ({\bf R}^k , Y).\label{NB-B}\end{equation}
(The isomorphism is not canonical but depends on the choice of a
linear basis $p_{\beta} (s)$ for the Weil algebra $\ci ({\bf R}^l )/I$.)

\section{Cahiers Topos}
    The site $\underline{D}$ of definition of this topos ${\cal C}$
is the dual of the category of $C^{\infty}$-rings of the form $\ci
({\bf R}^{l+k})/J$ where $J$ is a semi-Weil ideal, coming from a Weil
ideal $I\subseteq \ci ({\bf R}^l )$.  -- There is a full embedding $i:
\underline{Mf} \to {\cal C}$.

The full embedding $h$, described in \cite{CA}, of
$\underline{Con}^{\infty}$ into ${\cal C}$ is, on objects, given by
sending a CVS $X$ into the presheaf on $\underline{D}$ given by
$$\ci ({\bf R}^{l+k})/J \mapsto \ci ({\bf R}^{l+k}, X)/J(X).$$
For smooth maps $X\to Y$, composing with $Y$ preserves the property of
``being congruent mod $J$'', cf.\ Prop.\ 2.1 in \cite{CVSE} or Coroll.\  2 in
\cite{CA}, and this describes the functorality.  For finite
dimensional vector spaces $X$, $h(X)= i(X)$.

The embedding $h$ is full.  It preserves the exponentials in CVS, and
furthermore, if $X$ is a CVS, the $R$-module $h(X)$ in ${\cal C}$
``satisfies the vector form of Axiom 1'' (generalized Kock-Lawvere
Axiom), so that in particular synthetic calculus for curves $R\to
h(X)$ is available; cf.\ the final remark in \cite{CVSE}.  From this,
one may deduce that the embedding $h$ preserves differentiation, i.e.\
for $f:{\bf R}\to X$ a smooth curve, its derivative $f':{\bf R} \to X$
goes by $h$ to the synthetically defined derivative of the curve
$h(f): R=h({\bf R}) \to h(X)$.  This follows by repeating the argument
for Theorem 1 in \cite{PWAM} (the Theorem there deals with the case
where the codomain of $f$ is ${\bf R}$, but it is valid for $X$ as
well because $h(X)$ satisfies the vector form of Axiom 1).

\medskip

We note the following aspect of the embedding $h$.
Let $X$ be a CVS.  Each $\phi \in X'$ is  smooth linear $X\to {\bf
R}$ and hence defines a map $h(\phi ): h(X)\to h({\bf R}) = R$ in
${\cal C}$.  This map is $R$-linear.  \begin{prop}The maps $h(\phi ):
h(X)\to R$, as $\phi$ ranges over $X'$, form a jointly monic family.
\label{jointlymonic}\end{prop}
{\bf Proof.} The assertion can also be formulated: the natural map

$$e: h(X) \to \prod _{\phi \in X'} R$$
is monic (where $proj _{\phi} \circ e := h(\phi )$). To prove that
this (linear) map is monic, consider an element $a$ of the domain,
defined at stage $\ci ({\bf R}^{l+k})/J$, where $J$ is a semi-Weil
ideal. So $a\in \ci ({\bf R}^{l+k}, X)/J(X)$.
Let $\alpha \in  \ci ({\bf R}^{l+k}, X)$ be a smooth map representing
the class $a$, $a=\alpha + J(X)$.  The element $e(a)$ is the $X'$
tuple $a _{\phi} + J(X)$, where $a _{\phi} \in \ci ({\bf R}^{l+k})/J(X)$ is
represented by the smooth map $\phi \circ \alpha : {\bf R}^{l+k} \to
{\bf R}$.  To say $a$ maps to 0 by $e$ is thus to say that for each
$\phi \in X'$, $\phi \circ \alpha \in J$.  But this is precisely the
defining property for $\alpha$ itself to be in $J_w (X) = J(X)$, i.e.
for $a$ to be the zero as an element of $h(X)$ (at the given stage
$\ci ({\bf R}^{l+k})/J$).

\section{The internal space of test functions}

We first analyze the object $R^R$ in ${\cal C}$.  Because $h$
preserves exponentials, and $R= i({\bf R})=h({\bf R})$, $R^R$ is $h(\ci
({\bf R}))$.  Therefore an element of $R^R$ at stage $\ci ({\bf R}^{l+k})/J$,
where $J=p^* (I)$ is a semi-Weil ideal as above, is an element of
$$\ci ({\bf R}^{l+k},\ci ({\bf R}))/J(\ci ({\bf R})) \cong \prod
_B \ci ({\bf R}^k ,\ci ({\bf R}))\cong
\prod _B \ci ({\bf R}^{k+1}),$$ by (\ref{NB-B}).

In concrete terms, an element of $R^R$, defined at stage $\ci ({\bf
R}^{l+k})/J$, is thus given by the class mod $J^*$ of a smooth map $f:
{\bf R}^{l+k+1} \to {\bf R}$, $f(s,t,x)$, and even more concretely, by
the $B$-tuple of smooth maps ${\bf R}^{k+1}\to {\bf R}$, $D^{\beta} f
(0,t,x)$ (recall that the $D^{\beta}$'s differentiate in the
$s$-variable only, and then substitute $s=0$).

The following is a formula with a free variable $f$ that ranges over
$R^R$: \begin{equation} \exists b >0 [\; \; \forall x, (x<-b \vee x>b)
\Rightarrow f(x)=0\; \; ].
\label{old10}\end{equation}
Let us write $|x|>b$ as shorthand for the formula $x<-b \vee x>b$
    (so, in spite of the notation, we don't assume an ``absolute value''
    function).  Then the formula (\ref{old10}) gets the more readable
    appearance:
    \begin{equation} \exists b >0 [\; \; \forall x, |x|>b \Rightarrow
    f(x)=0\; \; ].  \label{fmla}\end{equation} (verbally: ``$f$ is a
    function $R\to R$ of bounded support'' (namely support contained in
    the interval $[-b,b]$).  Its extension is a subobject ${\cal
    D}(R)\subseteq R^R$.

We shall as a preliminary investigate when an element of $R^R$ defined
at a stage of the particular form $\ci ({\bf R}^k )$ belongs to the internal
${\cal D}(R)$ (described as the extension of the formula (\ref{fmla})).  So
for
the present, there are no Weil ideals invloved.

For simplicity of notation, let us write $K$ for $i({\bf R}^k ) =
R^k$.

So consider an element $f\in _K R^R$.  This means a map $K\to R^R$ in
${\cal C}$, and this in turn corresponds, by transposition, and by
fullness of the embedding $i$, to a smooth map
$$\hat{f} : {\bf R}^k \times {\bf R}\to {\bf R}.$$
Now we have that
$$\vdash _K \exists b >0 [\forall x, |x|>b \Rightarrow f(x)=0]$$
if and only if there is a covering $U_i$ of $K$ ($i\in I$) and
witnesses $b_i \in _{U_i} R_{>0}$, so that for each $i$
$$\vdash _{U_i} \forall x,\; |x|>b_i \Rightarrow f(x)=0$$
Externally, this implies that $b_i :U_i \to {\bf R}$ is a smooth
function with positive values, with the property that for all $t\in
U_i$, if $x$ has $x > b_i (t)$, then $\hat{f}(t,x)=0$.  The following
Lemma then implies that $f$ is of l.u.b.s.  on $U_i$, and since the
$U_i$'s cover $K$, $f$ is of l.u.b.s.  on $K$.

\begin{lemma} Let $g:U\times {\bf R} \to {\bf R}$ have the property
that there exists a smooth (or just continuous) $b: U\to R_{>0}$ so that
for all $t\in U$ $|x|>b(t)$ implies $g(t,x)=0$.  Then $g$ is l.u.b.s.
\label{lubs2}\end{lemma}
{\bf Proof.} For each $t\in U$, let $c_t$ denote $b(t)+1$.
There is a neighbourhood $V_t$ around $t$ such that $b(y) < c_t$ for
all $y\in V_t$.  The family of $V_t$'s, together with the constants
$c_t$ now witness that $g$ is l.u.b.s.  For, for all $y \in V_t$ and
any $x$ with $|x|>c_t$, we have $|x|> c_t > b(y)$, so $g(y,x)=0$.

Conversely, if $\hat{f}$ is l.u.b.s., it is easy to see that the
element $f\in _K R^R$ satisfies the formula (reduce to the uniformly 
bounded case, and write the condition as existence of a commutative 
square).

\medskip

So we conclude that for $f\in _K R^R$, $f\in _K {\cal D}(R)$ iff the
external function $f: K\times {\bf R}\to {\bf R}$ is l.u.b.s., i.e.,
by Theorem \ref{lubs}, iff $\hat{f}: K\to \ci ({\bf R})$ factors by a
(diffeologically!) smooth map through the inclusion ${\cal D}({\bf R})
\subseteq \ci ({\bf R})$, i.e.  belongs to $\ci ({\bf R}^k, {\cal
D}({\bf R})) = h( {\cal D}({\bf R}))(\ci ({\bf R}^k))$.  This proves
that, at least as far as generalized elements, defined at stages where
no Weil ideal is involved, we have ``$h({\cal D}({\bf R})) = {\cal
D}(R)$'', more precisely,

\begin{equation}h({\cal D}({\bf R}))(\ci ({\bf R}^k )) = {\cal
D}(R)(\ci ({\bf R}^k )).\label{compar}\end{equation}

To get a similar conclusion for elements of ${\cal D}(R)$ (as
synthetically defined by (\ref{fmla})),  defined at
stage $\ci ({\bf R}^{l+k})/J$, we shall prove that such can be represented
by $B$-tuples of elements defined at stage $\ci ( {\bf R}^k)$; we
shall prove that such a $B$-tuple defines an element of ${\cal D}(R)$
precisely if each of these $B$ elements is an element in ${\cal
D}(R)$.  This proof is a piece of ``purely synthetic reasoning'':

\medskip

    We consider an ${\bf R}$-algebra object $R$ in a topos ${\cal C}$,
and assume that $R$ satisfies the general ``Kock-Lawvere'' (K-L) axiom
(recalled
below), and is equipped with a strict order relation $<$.  Because the
reasoning is purely synthetic, we don't have to think in terms of
sheaf semantics, so for instance we don't have to be specific at what
``stages'', the ``elements'' in question are defined; we reason {\em
as if} all elements are global elements.  For $b>0$, we write $|x|>b$ as
shorthand for $x< -b \vee x>b$ as before; and we stress again that we
don't assume any absolute-value function (it does not exist in the
Cahiers topos).  We argue in ${\cal C}$ as if it were the category of
sets, making sure to use only intuitionistically valid reasoning.

A Weil algebra $\ci ({\bf R}^l )/I$, as above, gives rise to an
``infinitesimal'' subobject $W\subseteq R^l$: pick a (finite) set of
differential operators $D_{\beta}$ ($\beta \in B$) forming a basis for
$(\ci ({\bf R}^l )/I)^*$, and take the dual basis for $\ci ({\bf R}^l )/I$,
whose elements are represented mod $I$ by polynomials $p_{\beta} (s)$
in $l$ variables.  Then $W\subseteq R^l$ is the extension of the
formulas $p_{\beta} (s) =0$, $s$ being a variable ranging over $R^l$
(note that real polynomials in $l$ variables define functions $R^l \to
R$ in ${\cal C}$).

We assume that such $W$'s are internal atoms, in a sense we partially recall
below; this is so for all interesting models ${\cal C}$ of SDG,
including the Cahiers Topos.

To say that an $R$-module object $Y$ in ${\cal C}$ satisfies the general
K-L axiom is to say that for each such Weil algebra, the map
$$\prod _B Y \to Y^{W}$$
given by
$$(y_{\beta})_{\beta \in B} \mapsto [ s\mapsto \sum _B p_{\beta}(s)
\cdot y_{\beta}]$$
is an isomorphism.

We assume that $R$ itself satisfies K-L.  This immediately implies
that $R^M$ does for any $M\in {\cal C}$.  We shall consider $R^R$.

Now recall that
${\cal D}(R) \subseteq R^R$ was the subobject which is the extension of
the formula (\ref{fmla}) (with free variable $f$ ranging over $R^R$)
$\exists b >0: |x|>b \Rightarrow f(x)=0$.

\begin{prop} Let a $B$-tuple of elements $f_{\beta}$ in $R^R$ represent an
element in $(R^R)^W$.  Then it defines an element in the sub``set''
$({\cal D}(R))^W$ if and only if each $f_{\beta}$ is in ${\cal D}(R)$.
\label{synthetic}\end{prop}

{\bf Proof.} Assume first that all
$f_{\beta}$ are in ${\cal D}(R)$.  For each $\beta$ there exists a
witnessing $b_{\beta}>0$ witnessing that the formula (\ref{fmla})
holds for $f_{\beta}$, but since there are only finitely many
$\beta$'s, we may assume one common witness $b>0$.  So for all
$\beta$, and for all $x$ with $|x|>b$, $f_{\beta }(x)=0$.  But then
for each such $x$, the function of $s\in W$ given by
$$s\mapsto \sum _{\beta} p_{\beta }(s) \cdot f_{\beta}(x) $$
is the zero function.  The sum here, as a function of $s$ and $x$, is the
element
of $(R^R)^W$ corresponding to the $B$-tuple $f_{\beta}$, and for
$|x|>b$, it is the zero.  So for each $s$, the given fixed $b$ witnesses
that the sum, as a function of $x$, is in ${\cal D}(R)$.

Conversely, assume that the $f_{\beta}$'s are such that the
corresponding function $W\to R^R$ factors through ${\cal D}(R)$.  So
for each $s\in W$, the function
$$x\mapsto \sum _{\beta} p_{\beta }(s) \cdot f_{\beta}(x)$$
belongs to ${\cal D}(R)$.
So
\begin{equation}\forall s\in W \; \exists b >0 : |x|>b \Rightarrow
\sum _{\beta} p_{\beta }(s) \cdot f_{\beta}(x) =0.\label{Sko}\end{equation}
  We
would like to {\em pick} for each $s\in W$ a $\tilde{b}(s)$ such that
$$\forall s\in W: |x|>\tilde{b}(s)
\Rightarrow \sum _{\beta} p_{\beta }(s) \cdot f_{\beta}(x) =0;$$
the existence of such a {\em function} $\tilde{b}$ follows from
(\ref{Sko}) by a use of the Axiom of Choice, so in general is not
possible in a topos.  But since $W$ is an internal atom, and $s$
ranges over $W$, such a function $\tilde{b}$ exists after all.  (See
the Appendix for a general formulation and proof of this principle.)

But now $|x|>\tilde{b}(0) \Rightarrow |x| > \tilde{b}(s)$ for all
$s\in W$, because $\tilde{b}$, as does any function, preserves
infinitesimals, and because strict inequality is unaffected by
infinitesimals.  So we have a $b$, namely $\tilde{b}(0)$, so that
$$\forall s\in W: |x|>b \Rightarrow \sum p_{\beta}(s)\cdot f_{\beta }
(x) =0.$$ So for $|x|>b$, $$\forall s\in W , \sum p_{\beta}(s)\cdot
f_{\beta }
(x) =0.$$
Thus, for fixed $x$ with $|x|>b$, the function of $s$ here is 
constantly
0.  But functions $W\to R$ can uniquely be described as linear
combinations of the $p_{\beta}(s)$'s (this is a verbal rendering of
the K-L axiom for $R$).  So for such $x$ each $f_{\beta }(x)$ is $0$.
So $b$ witnesses, for each $\beta$, that $f_{\beta } \in {\cal D}(R)$.
This proves the Proposition.

\medskip

Combining (\ref{NB-B}) (with ${\cal D}({\bf R})$ for $Y$) with
(\ref{compar}) and Proposition \ref{synthetic}, we get

\begin{thm} The subobject ${\cal D}(R)$ of $R^R$ is exactly $h({\cal
D}({\bf R}))$.
\label{boundedsupp}\end{thm}

\medskip

The $R$-module ${\cal D}(R)\subseteq R^R$ (=the extension of the
formula (\ref{fmla}))  is by definition the internal vector space of test
functions; and we form the subobject $${\cal D}'(R) \subseteq
R^{{\cal D}(R)}$$
which is the extension of the formula ``$\phi$ is $R$-linear''
($\phi$ a variable ranging over $R^{({\cal D}(R)}$).  So ${\cal D}'(R)$ is the
internal vector space of distributions.

\medskip

We make an analysis of $h(Y')$ for a general CVS $Y$.  Recall that the
diffeology on $Y'$ is inherited from that of $\ci (Y,{\bf R})$, so
that (for an open $U\subseteq {\bf R}^k$), the smooth plots $U\to Y'$
are in bijective correspondence with smooth maps $U\times Y \to {\bf
R}$, which are ${\bf R}$-linear in the second variable $y\in Y$.  It
follows that the elements at stage $\ci ({\bf R}^k )$ (no Weil ideal
involved) are in bijective correspondence with smooth maps ${\bf R}^k
\times Y \to {\bf R}$, ${\bf R}$-linear in the second variable, or
equivalently, with smooth ${\bf R}$-linear maps $Y\to \ci ({\bf R}^k , {\bf
R})$, i.e.\ with smooth ${\bf R}$-linear maps $Y\to \ci ({\bf R}^k )$.

On the other hand, an element of $R^{h(Y)}$ defined at stage $\ci ({\bf
R}^k)$ is a morphism
$R^k \to R^{h(Y)}$, hence by double transposition it corresponds to a
map $h(Y)\to  R^{R^k }$; and it belongs to the subobject $Lin _R (h(Y), R)$
iff its double transpose is $R$-linear.  Since $h$ is full and
faithful, and preserves the cartesian closed structure (hence the
transpositions), this double transpose corresponds bijectively to a
smooth map $Y\to \ci ({\bf R}^k ,{\bf R})= \ci ({\bf R}^k)$, and
$R$-linearity is equivalent to ${\bf R}$-linearity, by the following general
\begin{lemma} Let $X$ and $Y$ be CVS's.  Then a smooth map $f:
Y\to X$ is ${\bf R}$-linear iff $h(f): h(Y) \to h(X)$ is $R$-linear.
\end{lemma}
{\bf Proof.} The implication $\Rightarrow$ is a consequence of the
fact that $h$ preserves binary cartesian products (and of $h({\bf
R})=R$).  For the implication $\Leftarrow$, we just apply the global
sections functor $\Gamma$; note that $\Gamma (Y)$ is the underlying
set of the vector space $Y$, and similar for $X$; and $\Gamma (R) =
{\bf R}$.

\medskip

    We have in particular:

\begin{prop}There is a natural one-to one correspondence between
distributions on ${\bf R}$, and $R$-linear maps ${\cal D}(R)\to R$
\label{com}\end{prop}
{\bf Proof.} By fullness of the embedding $h$ of
$\underline{Con}^{\infty}$ into the Cahiers topos ${\cal C}$, there is a
bijection between the set of smooth maps ${\cal D}({\bf R}) \to {\bf
R}$, and the set of morphisms in ${\cal C}$, $h ({\cal D}({\bf R}))
\to h({\bf R}) =R$, and ${\bf R}$-linearity corresponds to
$R$-linearity, by the above Lemma. The result now follows from
$h({\cal D}({\bf R}) )= {\cal D}(R)$ (Theorem \ref{boundedsupp}).

\medskip This result should be compared to the Theorem of
\cite{QR}, or Proposition II.3.6 in \cite{MR}, where a related
assertion is made for distributions-with-compact-support, i.e.\ where
${\cal D}(R)$ is replaced by the whole of $R^R$, -- or even with
$R^M$, with $M$ an arbitrary smooth manifold.  Distributions with
compact support are generally easier to deal with synthetically (as we
did in \cite{WE}), but they are not adequate for the heat equation.

\section{Half Line}

By Theorem \ref{p33}, the two $\ci$-rings $\ci ({\bf R})/{\cal
M}^{\infty}_{\geq 0}$ and $\ci ({\bf R}_{\geq 0})$ are isomorphic, where
${\cal
M}^{\infty}_{\geq 0}$ is the ideal of smooth functions vanishing on the
non-negative half line, and $\ci ({\bf R}_{\geq 0})$ is the ring of
smooth functions ${\bf R}_{\geq 0} \to {\bf R}$.  Being a quotient
of the ring $\ci ({\bf R})$ which represents $R \in {\cal C}$, it
defines a subobject of $R$, which we denote $R_{\geq 0}$ (also
considered in \cite{KR1}\footnote{The ring representing $R_{\geq 0}$,
   was in loc.cit.\ defined using the ideal ${\cal M}_{\geq 0}^g$ of
functions vanishing on an
open neighbourhood of ${\bf R}_{\geq 0}$, rather than the ideal ${\cal
M}_{\geq 0}^{\infty}$
considered here. But it can be proved that they represent (from the
outside) the same object in the Cahiers topos.     }).
-- Thus,  $R_{\geq 0}$ is ``represented
from the outside'' by the $\ci $-ring  $\ci ({\bf R})/{\cal
M}^{\infty}_{\geq 0} \cong\ci ({\bf R}_{\geq 0})$.

\begin{prop}
Let $I\subseteq C^{\infty}({\bf R}^l)$ be a Weil ideal and let $f:
{\bf R}^l\times {\bf R}^n\lra {\bf R}$ be a smooth function. Then the
following are equivalent:
\begin{enumerate}
\item  $f(\underline{0}, x)\geq 0$ for all $x\in {\bf R}^n$
\item $\rho (f(\underline{w},x))\in I^*$  for all $\rho \in
m_{{\bf R}\geq 0}^{\infty}.$
\end{enumerate}
\end{prop}
{\it Proof:} ``not 1'' implies ``not 2''; for, if $f(0,x)<0$, we may
find a function $\rho$ vanishing on ${\bf R}_{\geq 0}$ and with
value 1 at $f(0,x)$.  Then $f \notin I^*$ (recall that any Weil ideal
$I$ consists of functions vanishing at $0$).

\medskip
1 implies 2: By Taylor expansion,

$\begin{array}{lll} \\
(\rho\circ f)(\underline{w},x) & = & (\rho\circ f)(\underline{0},x) +
\sum_iw_i(\rho\circ f)'_i (\underline{0}. x) + \\
&  & \sum_{i,j}w_iw_j(\rho\circ f)'_{i,j}(\underline{0}, x) + \dots
\end{array}$

\medskip
where $(-)_i= \partial/\partial x_i,\; (-)_{i,j}=\partial^2/\partial
x_ix_j$ etc.

\medskip
This series finishes after finitely many terms modulo $I^*,$ since a
product of powers of $w_i$'s belong to the ideal $I$. But each of its
terms is 0: Indeed, so is the term without
derivatives, by hypothesis. But so are the others. For instance.
$(\rho\circ f)'_i (\underline{0}, x)=\rho'(f(\underline{0}, x)\partial
f/\partial x_i(\underline{0}, x)$ is 0, since the derivative of $\rho$
is zero on non-negative reals (by definition of $m_{{\bf R}\geq
0}^{\infty}$).

\medskip

An element $F$ of $R_{\geq 0}$ defined at stage $\ci ({\bf
R}^{l+k})/J$ is represented by a function $f$ satisfying the
conditions of the Proposition.

\begin{prop} There is a bijection between the set of maps
    $K: {\bf R}_{\geq 0} \to X$ which are strongly smooth 
    (in the sense of Definition \ref{strong} in Section 1),
    and the set of maps
$\overline{K}: R_{\geq 0} \to h(X)$ 
in ${\cal C}$.  \label{bar}\end{prop} {\bf
Proof/Construction.} Given an element $F$ of $R_{\geq 0}$ defined at
stage \newline $\ci ({\bf R}^{l+k})/J$, represented by $f$, as above.
We want to produce an element $K(F)$ of $X$ defined at stage $\ci
({\bf R}^{l+k})/J$, in other words, an element of $\ci ({\bf R}^{l+k},
X)/J(X)$.  We take the element represented (mod $J(X)$) by the smooth
map ${\bf R}^{l+k} \to
X$ given by the $r$ first terms of the series $(s,t)\mapsto$ $$
K(f(0,t)) + K'(f(0,t))\cdot [f(s,t)-f(0,t)] +
\frac{K''(f(0,t))}{2!} \cdot [f(s,t)-f(0,t)]^2 + \ldots .$$
Note that $[f(s,t)-f(0,t)]^r \in J$ since for fixed $t$,
$f(s,t)-f(0,t) \in {\cal M}$, hence $[f(s,t)-f(0,t)]^r \in {\cal M}^r
\subseteq I$.  This suffices, by the description of semi-Weil ideals
in terms of differential operators.

We have to prove that the class mod $J(X)$ of this map only depends on
the class of $f$.  If we change $f$ into $f+h$ with $h\in J$, the
term in the square brackets (real numbers!) change into $[f(s,t)-f(0,t) +
h(s,t)]$ (using that $h(0,t)=0$), and the ``Taylor coefficients''
($\in X$) do not change at all, for the same reason.  -- Uniqueness is
easy, using Proposition \ref{jointlymonic}, together with the fullness 
result from \cite{PR} on manifolds with boundary.

\medskip The Proposition is a ``mixed fullness'' result; we have that
   $Con ^{\infty}$ and $\underline{Mf}$ (= smooth manifolds), (even
   the category of smooth manifolds with boundary), embed fully in The Cahiers
   Topos; but at present we do not  have a general result about what can
   be said about $\ci (M, X)$, for $M$ a manifold (possibly with
   boundary) and $X$ a CVS -- not to speak of $\ci (X, M)$.

      \section{Heat Equation in the Cahiers Topos}

      For any topos ${\cal C}$ with a ring object $R$ with a preorder
      $\leq$, we may form the $R$-module ${\cal D}'(R^n )$ of distributions on
      $R^n$, as explained in Section 7.  If ${\cal C}, R$ is a model of
      SDG, then ${\cal D}'(R^n )$ automatically satisfies the ``vector
      form'' of the general Kock-Lawvere axiom, so that (synthetic)
      differentiation of functions $K: R \to {\cal D}'(R^n )$ is possible -
      it is even enough that $K$ be defined on suitable
      (``formally etale'') subobjects of $R$, like $R_{\geq 0}$.  We think of
the
      domain $R$ or $R_{\geq 0}$ as ``time'', and denote the differentiation
      of curves $K$ w.r.\ to time by the Newton dot, $\dot{K}$.  On the
      other hand, we think of $R^n$ as a space, and the various partial
      derivatives $\partial /\partial x_i$ ($i=1, \ldots ,n$), as well as
      their iterates, we call spatial derivatives; in case $n=1$, they are
      just denoted $(-)'$, $(-)''$, etc.  They live on ${\cal D}'(R^n )$ as
      well, by the standard way of differentiating distributions (which
      immediately translates into the synthetic context, cf.\ e.g.\
      \cite{WE}).  The heat equation for (Euclidean) space in $n$ dimensions
      says $\dot{K} = \Delta \circ K$, where $\Delta$ is the Laplace
      operator; in one dimension it is thus the equation
      $$\dot{K} = K''.$$

      We can summarize the constructions into an general existence theorem
      about models for SDG:

      \begin{thm}There exists a well-adapted model for SDG (with a preorder
      $\leq$ on $R$), in which the heat equation on the (unlimited) line $R$
      has a unique solution $k: R_{\geq 0} \to {\cal D}' (R)$ with initial value
$k(0) =
      \delta (0)$ (the Dirac distribution).
      \end{thm}
      {\bf Proof.} The well adapted model witnessing the validity of the
      Theorem is the Cahiers Topos ${\cal C}$. Consider the classical heat
kernel,
      viewed, as we did in Section 4, as a map ${\bf R}_{\geq 0} \to {\cal
      D}' ({\bf R})$.  By Theorem \ref{strongsmooth}, this map is smooth in the
      strong sense, hence by Proposition \ref{bar}, it
      defines a morphism in ${\cal C}$, $\bar{K}: R_{\geq 0} \to h({\cal
      D}'({\bf R})$.  This $\bar{K}$ is going to be our $k$.  By Theorem
      (\ref{boundedsupp})), its codomain is the desired ${\cal D}'(R)$.  So
      all that remains is to prove that this $k$ satisfies the heat equation
      $\dot{k} = \Delta \circ k$.  This is a purely formal argument from the
      fact that $K$ does, and the fact that $h$ takes ``analytic''
      differentiation into the ``synthetic'' differentiation in ${\cal C}$.
      We give this argument.  Synthetically, we want to prove that for all
      $x\in {\bf R}_{\geq 0}$ and $d\in D$
      $$k(x+d) = k(x) + d\cdot  \Delta (k(x)).$$
    Universal validity of this equation means that a certain diagram, with
    domain $R_{\geq 0} \times D$ and codomain ${\cal D}'(R)$, commutes.
    Taking the transpose of this diagram, we get a diagram with domain
    $R_{\geq 0}$ and codomain $({\cal D}'(R))^D \cong {\cal D}'(R)\times
    {\cal D}'(R)$ (by K-L for ${\cal D}'(R)$):
    \begin{diagram}R_{\geq 0} & \rTo ^{\hat{+}} & (R_{\geq 0})^D\\
      \dTo^k && \dTo _{k^D}\\
      {\cal D}'(R)&&&\\
      \dTo^{(1,\Delta )}&&&\\
      {\cal D}'(R)\times {\cal D}'(R)&\lTo_{\cong} &({\cal D}'(R))^D
      \end{diagram}

      When the global sections functor $\Gamma$ is applied to this diagram,
      the left hand column yields $(K,\Delta \circ K)$, because $\Gamma (k)=K$;
      the composite of the other maps is $(K, \dot{K})$ because $\Gamma$
      takes synthetic differentiation into usual differentiation.  Since $K$
      satisfies $\dot{K}= \Delta \circ K$, we conclude that $\Gamma$ applied
      to the exhibited diagram commutes.  Now $\Gamma$ is not faithful, but
      because of the special form of the domain and codomain of the two maps
      to be compared, we may still get the conclusion, by virtue of the
following

    \begin{prop} Given a map $a: R_{\geq 0}\to h(X)$, where $X$ is a CVS.
    If $\Gamma (a) =0$, then $a=0$.
    \end{prop}
    {\bf Proof.} Since the $h(\phi ): h(X)\to R$ are jointly monic as
    $\phi$ ranges over $X'$, by Proposition \ref{jointlymonic}, it suffices to
    see that each $h(\phi )\circ a$ is $0$.  Since $\Gamma (h(\phi )\circ
    a)= \phi \circ \Gamma (a)$, this reduces the question to the case
    where $X={\bf R}$.  A map $a: R_{\geq 0} \to R$ is tantamount to an
    element in $\bar{a}: \ci ({\bf R}_{\geq 0})$, and the assumption $\Gamma
    (a)=0$ is tantamount to $\bar{a} (t)=0$ for all $t\in {\bf R}_{\geq
    0}$.  But this clearly implies that $\bar{a}$, and hence $a$, is $0$.

\section*{Appendix}
Recall that an {\em atom} $A$ in a cartesian closed category ${\cal C}$ is an
object so that the exponential functor $(-)A$ has a right adjoint; in
particular, it takes epimorphisms to epimorphisms.  The following says
that ``axiom of choice'' holds for ``$A$''-tuples  sets:
\begin{prop}
Assume that $A$ is an atom, $B$ arbitrary and $R\subseteq A\times B.$
Then
$$(\forall a\in A)(\exists b\in B)\; R(a,b)\Longrightarrow (\exists
\tilde{b}\in B^A) (\forall a\in A)\;R(a, \tilde{b}(a))
$$
\end{prop}
{\it Proof: } The hypothesis means that the composite $R\to A\times
B\stackrel{\pi_1}\lra A$ is surjective.  By exponentiation, and the
assumption that $A$ is an atom,  the composite $R^A\to A^A\times
B^A\stackrel{\pi_1}\lra A^A$ is surjective. In particular, $1_A\in A^A$
must have a pre-image $(1_A,\tilde{b}).$ This $\tilde{b}$ obviously
does the job.

\end{document}